\documentclass{article}
\usepackage{amsmath, amssymb, amsthm}
\usepackage{caption} 
\usepackage{hyperref}
\usepackage{setspace} 
\usepackage{tikz}
\usetikzlibrary{decorations.pathreplacing} 

\usepackage[
  block=ragged
]{biblatex}
\newtheorem{theorem}{Theorem}
\newtheorem{conjecture}[theorem]{Conjecture}
\newtheorem{definition}[theorem]{Definition}
\newtheorem{example}[theorem]{Example}
\newtheorem{note}[theorem]{Note}
\newtheorem{openquestion}[theorem]{Open Question}
\newtheorem{proposition}[theorem]{Proposition}

\begin{document}

\title{Spinning switches on a wreath product}
\author{Peter Kagey}

\maketitle

\begin{abstract}
In this paper, we attempt to classify an algebraic phenomenon on wreath products
that can be seen as coming from a family of puzzles about switches on the corners
of a spinning table. Such puzzles have been written about and generalized since they
were first popularized by Martin Gardner in 1979.
In this paper, we provide perhaps the fullest generalization yet, modeling
both the switches and the spinning table as arbitrary finite groups combined
via a wreath product. We classify large families of wreath products depending on
whether or not they correspond to a solvable puzzle,
completely classifying the puzzle in the case when the switches behave like abelian groups,
constructing winning strategies for all wreath product that are $p$-groups,
and providing novel examples for other puzzles where the switches behave like nonabelian groups,
including the puzzle consisting of two interchangeable copies of the monster group $M$.
Lastly, we provide a number of open questions and conjectures, and
provide other suggestions of how to generalize some of these ideas further.
\end{abstract}

\section{Overview and preliminaries}
\label{sec:overviewAndPreliminaries}
This paper is organized into six sections.
This section,
Section \ref{sec:overviewAndPreliminaries} provides a brief history of this genre of puzzles and introduces some of the first approaches to generalizing the puzzle further.
Section \ref{sec:WreathModel} models these generalizations in the context of the wreath product, and formalizes the notation of puzzles being solvable.
Section \ref{sec:Reductions} explores situations where the puzzle does not have a winning strategy, and provides reductions that allow us to prove that entire families of puzzles are not solvable.
Section \ref{sec:pGroupStrategy} constructs a strategy for switches that behave like $p$-groups, and gives us ways of building strategies from smaller parts.
Section \ref{sec:OtherSurjectiveStrategies} provides novel examples of puzzles that do not behave like $p$-groups, but still have winning strategies.
Lastly,
Section \ref{sec:OpenQuestions} provides further generalizations, and contains dozens of conjectures, open questions, and further directions.

\subsection{History}
Generalized spinning switches puzzles are a family of closely related puzzles
that were first popularized by Martin Gardner  in a puzzle called
``The Rotating Table''
in the February 1979 edition of his column ``Mathematical Games'' \cite{Gardner1979Problem}.
Gardner writes that he learned of the puzzle from Robert Tappay of Toronto who
``believes it comes from the U.S.S.R.,'' a history that is not especially
forthcoming.

My preferred version of the puzzle appears in Peter Winkler's 2004 book
\textit{Mathematical Puzzles A Connoisseur's Collection} \cite{Winkler2004}
\begin{quote}
  Four identical, unlabeled switches are wired in series to a light bulb.
  The switches are simple buttons whose state cannot be directly observed,
  but can be changed by pushing; they are mounted on the corners of a
  rotatable square. At any point, you may push, simultaneously, any subset
  of the buttons, but then an adversary spins the square. Show that there
  is a deterministic algorithm that will enable you to turn on the bulb in
  at most some fixed number of steps.
\end{quote}

(Winkler's version will be a working example in many parts of this paper, so it
is worth keeping in mind.
An illustration can be found in Figure \ref{fig:twoSwitches}.)

Over the last three decades, various authors have considered generalizations of
this puzzle. Here, we build on those results and go further.
The first place authors looked to generalize was suggested by Gardner himself.
In his March 1979 column the following month, he provided the answer to the
original puzzle and wrote
\begin{quote}
The problem can also be generalized by replacing glasses with objects that
have more than two positions. Hence the rotating table leads into deep
combinatorial questions that as far as I know have not yet been explored.
\cite{Gardner1979Solution}
\end{quote}

In 1993, Bar Yehuda, Etzion, and Moran \cite{BarYehuda1993}. took on the challenge and developed a theory
of the spinning switches puzzle where the switches behave like roulettes
with a single ``on'' state. In this paper we take Gardner's charge to
its logical conclusion and consider switches that behave like arbitrary
``objects that have more than two positions''.

Another generalization of this puzzle could look at other ways of ``spinning''
the switches. In 1995, Ehrenborg and Skinner \cite{Ehrenborg1995} did this
in a puzzle they call "Blind Bartender with Boxing Gloves", that analyzed
this puzzle while allowing the adversary to use an arbitrary, faithful group action
to ``scramble'' the switches. We analyze our generalized switches within this
same context.

This puzzle was re-popularized in 2019 when it appeared in ``The Riddler''
column from the publication FiveThirtyEight \cite{FiveThirtyEight}.
Shortly after this, in 2022, Yuri Rabinovich synthesized Bar Yehuda and Ehrenborg's
results in a paper that modeled the collection of switches as a vector space
over a finite field, and modeled the ``spinning'' or ``scrambling'' as a
faithful, linear group action on this vector space.

For more background, see Sidana's \cite{Sidana2020} detailed overview of the
history of this and related problems.

\subsection{A solution to Winkler's Spinning Switches puzzle}

We will start by discussing the solution to Winkler's version of the puzzle
because the solution provides some insights and intuition for the techniques
that we use later. Before solving the four-switch version of the puzzle,
we will make Peter Winkler proud by beginning with a simpler, two-switch version.

\begin{figure}
  \center
  \begin{tikzpicture}[scale=0.8]
    \draw (0,0) rectangle (4,4);
    \draw (2.1,1.9)--(1.9,2.1);
    \draw (2.1,2.1)--(1.9,1.9);
    \draw (1,3) circle (0.5);
    \draw (3,1) circle (0.5);
    \draw[ultra thick, ->] (0.3, 4.3) to[out=45, looseness=1.7, in=45] (4.3,0.3);
    \draw[ultra thick, ->] (3.7,-0.3) to[out=225, looseness=1.7, in=225] (-0.3, 3.7);

    \draw[xshift=-18em] (0,0) rectangle (4,4);
    \draw[xshift=-18em] (2.1,1.9)--(1.9,2.1);
    \draw[xshift=-18em] (2.1,2.1)--(1.9,1.9);
    \draw[xshift=-18em] (1,1) circle (0.5);
    \draw[xshift=-18em] (1,3) circle (0.5);
    \draw[xshift=-18em] (3,1) circle (0.5);
    \draw[xshift=-18em] (3,3) circle (0.5);
    \draw[xshift=-18em, ultra thick, ->] (0.3, 4.3) to[out=45, looseness=0.9, in=135] (3.7,4.3);
    \draw[xshift=-18em, ultra thick, ->] (4.3,3.7) to[out=-45, looseness=0.9, in=45] (4.3,0.3);
    \draw[xshift=-18em, ultra thick, ->] (3.7, -0.3) to[out=-135, looseness=0.9, in=-45] (0.3,-0.3);
    \draw[xshift=-18em, ultra thick, ->] (-0.3, 0.3) to[out=-225, looseness=0.9, in=-135] (-0.3,3.7);
  \end{tikzpicture}
  \caption[Illustration of Winkler's Spinning Switches and a two-switch analog.]{
    An illustration of Winkler's Spinning Switches puzzle and a two-switch analog.
  }
  \label{fig:twoSwitches}
\end{figure}
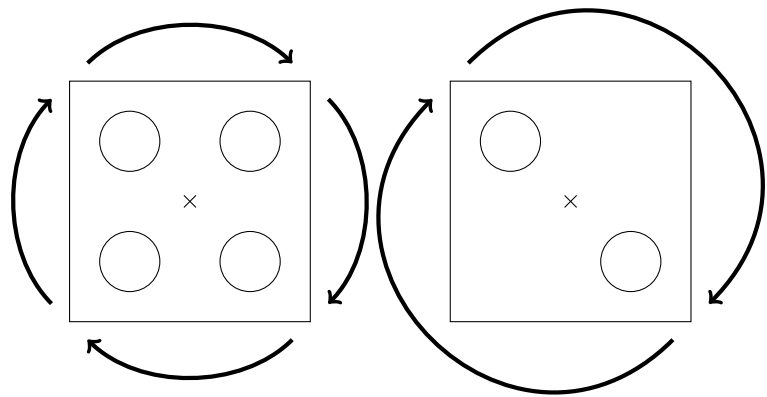

\begin{example}
  Suppose that we have two identical unlabeled switches on opposite corners
  of a square table, as in Figure \ref{fig:twoSwitches}

  Then we have a three-step solution for solving the problem.
  We start by toggling both switches simultaneously, and allow the adversary
  to spin the table.
  If this does not turn on the light,
  this means that the switches were (and still) are in different states.

  Next, we toggle one of the two switches
  to ensure that the switches are both in the same state. If the light has not
  turned on, both must be in the off state.

  The adversary spins the table once more, but to no avail. We know both
  switches are in the off state, so we toggle them both simultaneously, turning
  on the lightbulb.
\end{example}

In order to bootstrap the two-switch solution into a four-switch solution,
we must notice two things: \begin{enumerate}
  \item First, if we can get two switches along each diagonal into the same state
  respectively, then we can solve the puzzle by toggling both diagonals
  (all four switches),
  followed by both switches in a single diagonal,
  and lastly both diagonals again.
  In this (sub-)strategy, toggling both switches along a diagonal is
  equivalent to toggling a single switch in the two-switch analog.

  \item Second, we can get both diagonals into the same state at some point
  by toggling a switch from each diagonal (two switches on any side of the square),
  followed by a single switch from one diagonal,
  followed by again toggling a switch from each diagonal.
\end{enumerate}

We will interleave these strategies in a particular way, following the notation
of Rabinovich \cite{Rabinovich2022}.

\begin{definition}
  Given two sequences $A = \{a_i\}_{i=1}^N$ and $B = \{b_i\}_{i=1}^M$, we can
  define the \textbf{interleave} operation as \begin{align}
    A \circledast B &= (A,b_1,A,b_2,A,\dots,b_M,A) \\
      &= (
      \underbrace{a_1, \dots, a_N}_A,
      b_1,
      \underbrace{a_1, \dots, a_N}_A,
      b_2,
      \underbrace{a_1, \dots, a_N}_A,
      \dots,
      b_M,
      \underbrace{a_1, \dots, a_N}_A).
  \end{align} which has length $(M+1)N + M = MN + M + N$.
\end{definition}

Typically it is useful to interleave two strategies when
$A$ solves the puzzle given that the switches are in a particular state, and
$B$ gets the switches into that particular state.
We also need $A$ not to ``interrupt'' what $B$ is doing.
In the problem of four switches on a square table,
$B$ will ensure that the switches are in the same state within each diagonal,
and $A$ will turn on the light when that is the case.
Moreover, $A$ does not change the state \textit{within} either diagonal.

\begin{proposition}
  There exists a fifteen-move strategy that guarantees that the light in
  Winkler's puzzle turns on.
  \label{prop:WinklersSolution}
\end{proposition}
\begin{proof}
  We begin by formalizing the two strategies. We will say that the first
  strategy $S_1$ where we toggle the two switches in a diagonal together
  will consist of the following three moves: \begin{enumerate}
    \item Switch \textbf{a}ll of the bulbs ($A$).
    \item Switch the \textbf{d}iagonal consisting of the upper-left and lower-right bulbs ($D$).
    \item Switch \textbf{a}ll of the bulbs ($A$).
  \end{enumerate}
  We will say that the second strategy $S_2$ where we get the two switches
  within each diagonal into the same state consists of the following three
  moves: \begin{enumerate}
    \item Switch both switches on the left \textbf{s}ide ($S$).
    \item Switch \textbf{one} switch ($1$).
    \item Switch both switches on the left \textbf{s}ide ($S$).
  \end{enumerate}
  Then the $15$ move strategy is \begin{equation}
    S_1 \circledast S_2 = (A, D, A, S, A, D, A, 1, A, D, A, S, A, D, A)
  \end{equation}
\end{proof}

We will generalize this construction in
Theorem \ref{thm:surjectiveStrategyDecomposition},
which offers a formal proof that this strategy works.

(It is worth briefly noting that $S_1 \circledast S_2$ is the fourth
\textit{Zimin word} (also called a \textit{sequipower}),
an idea that comes up in the study of combinatorics on words.)

\subsection{Generalizing switches}
\label{sec:GeneralizingSwitches}
Two kinds of switches are considered by Bar Yehuda, Etzion, and Moran in 1993
\cite{BarYehuda1993}: switches with a single ``on'' position that behave like
$n$-state roulettes ($\mathbb Z_n$) and switches that behave like
the finite field $\mathbb F_q$, both on a rotating $k$-gonal table.
We generalize this notion further by considering switches that behave like
arbitrary finite groups.

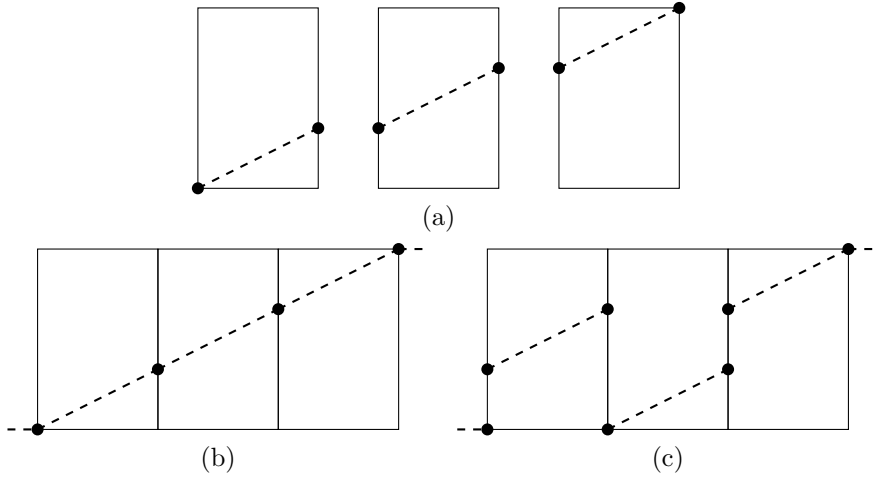
\begin{figure}
  \center{
  \begin{tikzpicture}[scale=0.8]
    \draw (0,0) rectangle ++(2,3);
    \draw[dashed,thick] (0,0) -- ++(2,1);
    \fill (0,0) circle (0.1);
    \fill (2,1) circle (0.1);

    \draw (3,0) rectangle ++(2,3);
    \draw[dashed,thick] (3,1) -- ++(2,1);
    \fill (3,1) circle (0.1);
    \fill (5,2) circle (0.1);

    \draw (6,0) rectangle ++(2,3);
    \draw[dashed,thick] (6,2) -- ++(2,1);
    \fill (6,2) circle (0.1);
    \fill (8,3) circle (0.1);
    \node at (4, -1/2) {(a)};
  \end{tikzpicture}
  }

  \noindent
  \begin{tikzpicture}[scale=0.8]
    \draw[dashed,thick] (-0.5,0) -- (0,0);
    \fill (0,0) circle (0.1);

    \draw (0,0) rectangle ++(2,3);
    \draw[dashed,thick] (0,0) -- ++(2,1);
    \fill (2,1) circle (0.1);

    \draw (2,0) rectangle ++(2,3);
    \draw[dashed,thick] (2,1) -- ++(2,1);
    \fill (4,2) circle (0.1);

    \draw (4,0) rectangle ++(2,3);
    \draw[dashed,thick] (4,2) -- ++(2,1);
    \fill (6,3) circle (0.1);

    \draw[dashed,thick] (6,3) -- (6.5,3);
    \node at (3, -1/2) {(b)};
  \end{tikzpicture}
  ~
  \begin{tikzpicture}[scale=0.8]
    \draw[dashed,thick] (-0.5,0) -- (0,0);
    \fill (0,0) circle (0.1);

    \draw (0,0) rectangle ++(2,3);
    \draw[dashed,thick] (0,1) -- ++(2,1);
    \fill (0,1) circle (0.1);
    \fill (2,2) circle (0.1);

    \draw (2,0) rectangle ++(2,3);
    \draw[dashed,thick] (2,0) -- ++(2,1);
    \fill (2,0) circle (0.1);
    \fill (4,1) circle (0.1);

    \draw (4,0) rectangle ++(2,3);
    \draw[dashed,thick] (4,2) -- ++(2,1);
    \fill (4,2) circle (0.1);
    \fill (6,3) circle (0.1);

    \draw[dashed,thick] (6,3) -- (6.5,3);
    \node at (3, -1/2) {(c)};
  \end{tikzpicture}
  \caption[A schematic for a switch that behaves like $S_3$.]{
    Part (a) shows a simple schematic for the components of a switch that
    behaves like $S_3$, the symmetric group on three letters.
    The three rectangles can be permuted arbitrarily, but only configuration (b)
    completes the circuit. All other configurations fail to
    complete the circuit (e.g., (c)).
  }
  \label{fig:S3Switch}
\end{figure}

\begin{example}
In Figure \ref{fig:S3Switch}, we provide a schematic for a switch that behaves
like the symmetric group $S_3$.
It consists of three identical-looking parts that need to be
arranged in a particular order in order for the switch to be on.

We could also construct a switch that behaves like the dihedral group of the
square, $D_8$.
This switch might look like a flat, square prism that can slot into a square hole,
such that only one of the $|D_8| = 8$ rotations of the prism completes the circuit.
\label{ex:S3D8Schematics}
\end{example}

\begin{note}
  One subtlety of using a group $G$ to model a switch is that
  both the ``internal state'' of a switch itself and
  the set of ``moves'' or changes are modeled by $G$.
  Therefore it may be useful to think of the state as the underlying set of $G$
  where the moves act via a right group action of $G$ on itself.
\end{note}

The reason that it is appropriate to use a group to model a switch
is because groups have many of the properties we would expect in a desirable
switch.
\begin{note}
  The axioms for a group $(G, \cdot)$ closely follow what we would expect from
  a switch.
\end{note}
\begin{enumerate}
  \item (Closure) The group $(G, \cdot)$ is equipped with a binary operation,
  $\cdot \colon G \times G \rightarrow G$. That is, for all pairs of elements
   $g_1, g_2 \in G$ their product is in $G$ \begin{equation}
     g_1 \cdot g_2 \in G.
  \end{equation}
  In the context of switches, this
  means that if the switch is in some state $g_1 \in G$ and the puzzle-solver
  applies the move $g_2 \in G$ to it,
  then the resulting state $g_1 \cdot g_2 \in G$ is in the set of possible
  states.
  \item (Identity) There exists an element $\operatorname{id}_G \in G$ such that
  for all $g \in G$, \begin{equation}
    \operatorname{id}_G \cdot g = g \cdot \operatorname{id}_G = g.
  \end{equation}
  This axiom is useful because it means that the puzzle-solver can ``do nothing''
  to a switch and leave it in whatever state it is in.
  Because the identity is a distinguished element in $G$,
  we will also use the convention that
  $\operatorname{id}_G$ is the ``on'' or ``winning'' state for a given switch.
  (It is worth noting that all of the arguments work basically the same way
  regardless of which element is designated as the on state.)
  \item (Inverses) For each element $g \in G$ there exists an inverse element
  $g^{-1} \in G$ such that \begin{equation}
    g \cdot g^{-1} = g^{-1} \cdot g = \operatorname{id}_G.
  \end{equation}
  This axiom states that no matter what state a switch is in,
  there is a move that will transition it into the on state.
  \item (Associativity) Given three elements $g_1, g_2, g_3 \in G$,
  \begin{equation}
    (g_1 \cdot g_2) \cdot g_3 = g_1 \cdot (g_2 \cdot g_3)
  \end{equation}
  This axiom is not strictly necessary for modeling switches,
  but as we will see in a later definition, it gives us a convenient way to
  describe the conditions for a winning strategy.
  (In Subsection \ref{sub:quasigroupSwitches}, we briefly discuss dropping
  the associativity axiom by considering switches that behave like
  quasigroups with identity.)
\end{enumerate}


\subsection{Generalizing spinning}
We can also consider generalizations of ``spinning'' the switches.
In particular, we will adopt the generalization from
Ehrenborg and Skinner's \cite{Ehrenborg1995} 1995 paper, which use
arbitrary faithful group actions to permute the switches.
In particular, they provide a criterion that determines which group actions
yield a winning strategy in the case of a given number of ``ordinary'' switches
(those that behave like $\mathbb Z_2$).
Rabinovich \cite{Rabinovich2022} stretches these results a bit further and
looks at faithful linear group actions on collections of switches that
are modeled as a finite-dimensional vector space over a finite field.
We build on this result in the context of more general switches.

\section{A wreath product model}
\label{sec:WreathModel}
Recall that Peter Winkler's Spinning Switches puzzle consists of four two-way
switches on the corners of a rotating square table.
The behavior of the switches are naturally modeled as $\mathbb Z_2$, and
the rotating table is modeled as the cyclic group $C_4$.
The abstraction that takes these two groups and creates a model for the
underlying puzzle is the wreath product: that is, Winkler's puzzle behaves like
the wreath product of $\mathbb Z_2$ by $C_4$.

\subsection{Modeling generalized spinning switches puzzles}
We do not evoke wreath products arbitrarily: we use them because they are the
right abstraction to model a generalized spinning switches puzzle where
$G$ describes the behavior of the switches,
$\Omega$ describes the positions of the switches, and
the action of $H$ on $\Omega$ models the ways the adversary can permute the switches.

\begin{definition}[\cite{Rotman1999}]
  Let $G$ and $H$ be groups,
  let $\Omega$ be a finite $H$-set, and
  let $K = \prod_{\omega \in \Omega} G_\omega$, where $G_\omega \cong G$
  for all $\omega \in \Omega$.
  Then the \textbf{wreath product} of $G$ by $H$ denoted by $G \wr H$,
  is the semidirect product of $K$ by $H$,
  where $H$ acts on $K$ by $h \cdot (g_\omega) = g_{h^{-1}\omega}$ for $h \in H$
  and $g_\omega \in G_\omega$.
  The normal subgroup $K$ of $G \wr H$ is called
  the \textbf{base} of the wreath product.

  The group operation is $(k, h) \cdot (k', h') = (k(h \cdot k'), hh')$
\end{definition}

An element $(k, h) \in G \wr H$ represents a turn of the game:
The puzzle-solver chooses an element of the base $k \in K$ to indicate
how they want to modify each of their switches
and then the adversary chooses $h \in H$ and acts with $h$ on $\Omega$ to
permute the switches.

\begin{example}
  Consider the setup in the Winkler's Spinning Switches the puzzle, which
  consists of two-way switches ($\mathbb Z_2$) on the corners of a rotating
  square,
  $C_4 \cong \langle 0^\circ, 90^\circ, 180^\circ, 270^\circ \rangle$.
  The game itself corresponds to the wreath product $\mathbb Z_2 \wr C_4$.
  We will use the convention that the base of the wreath product, $K$, is
  ordered upper-left, upper-right, lower-right, lower-left, and that the group
  action is specified by degrees in the clockwise direction.

  Consider the following two turns:
  \begin{enumerate}
    \item During the first turn,
    the puzzle-solver toggles the upper-left and lower-right switches, and
    the adversary rotates the table $90^\circ$ clockwise.
    This is represented by the element \begin{equation}
      ((1,0,1,0), 90^\circ) \in \mathbb Z_2 \wr C_4.
    \end{equation}
    \item During the second turn,
    the puzzle-solver toggles the upper-left switch, and
    the adversary rotates the table $90^\circ$ clockwise.
    This is represented by the element \begin{equation}
      ((1,0,0,0), 180^\circ) \in \mathbb Z_2 \wr C_4.
    \end{equation}
  \end{enumerate}
  As illustrated in Figure \ref{fig:WreathProduct},
  the net result of these two turns is the same as
  a single turn where the puzzle-solver toggles the
  upper-left, upper-right, and lower-left
  switches and the adversary rotates the board $270^\circ$ clockwise.

  \begin{figure}
    \center
    \begin{tikzpicture}
      \draw[dashed] (-0.2,-0.2) rectangle (3.7,2.2);
      \draw [thick,decorate,decoration={brace,amplitude=10pt},yshift=-0.4pt]
        (3.7,-0.2) -- (-0.2,-0.2) node[black,midway,yshift=-0.6cm] {\footnotesize $((1,0,1,0), 90^\circ)$};

      \draw[dashed] (3.8,-0.2) rectangle (7.7,2.2);
      \draw [thick,decorate,decoration={brace,amplitude=10pt},yshift=-0.4pt]
        (7.7,-0.2) -- (3.8,-0.2) node[black,midway,yshift=-0.6cm] {\footnotesize $((1,0,0,0), 180^\circ)$};

      \draw[dashed] (7.8,-0.2) rectangle (10.2,2.2);
      \draw [thick,decorate,decoration={brace,amplitude=10pt},yshift=-0.4pt]
        (10.2,-0.2) -- (7.8,-0.2) node[black,midway,yshift=-0.6cm] {\footnotesize $((1,0,1,1), 270^\circ)$};

      \draw (0,0) rectangle (2,2);
      \node[circle,fill,text=white] at (1/2,3/2) {$s_1$};
      \node[circle,draw] at (3/2,3/2) {$s_2$};
      \node[circle,fill,text=white] at (3/2,1/2) {$s_3$};
      \node[circle,draw] at (1/2,1/2) {$s_4$};
      \draw[->] (2.5,1) -- node[above] {$90^\circ \! \curvearrowright$} (3.5,1) ;

      \draw (4,0) rectangle (6,2);
      \node[circle,fill,text=white] at (9/2,3/2) {$s_4$};
      \node[circle,draw,fill=gray] at (11/2,3/2) {$s_1$};
      \node[circle,draw] at (11/2,1/2) {$s_2$};
      \node[circle,draw,fill=gray] at (9/2,1/2) {$s_3$};

      \draw[->] (6.5,1) -- node[above] {$180^\circ \! \curvearrowright$} (7.5,1) ;

      \draw (8,0) rectangle (10,2);
      \node[circle,draw] at (17/2,3/2) {$s_2$};
      \node[circle,draw,fill=gray] at (19/2,3/2) {$s_3$};
      \node[circle,draw,fill=gray] at (19/2,1/2) {$s_4$};
      \node[circle,draw,fill=gray] at (17/2,1/2) {$s_1$};
    \end{tikzpicture}
    \caption[A wreath product interpretation of a spinning switches puzzle.]{
      An illustration of two turns in Winkler's Spinning Switches puzzle,
      modeled as elements of a wreath product.
    }
    \label{fig:WreathProduct}
  \end{figure}
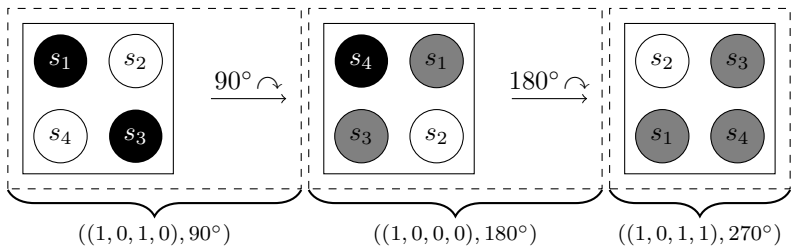

  The multiplication under the wreath product agrees with this: \begin{align*}
    ((1,0,1,0), 90^\circ) \cdot ((1,0,0,0), 180^\circ)
    &= ((1,0,1,0) + \underbrace{90^\circ \cdot (1,0,0,0)}_{(0,0,0,1)}, 90^\circ + 180^\circ) \\
    &= ((1,0,1,1), 270^\circ)
  \end{align*}
\end{example}

As suggested earlier, it is occasionally useful to designate a particular state
of the switches as the winning state. We will use the
convention that the lightbulb turns on when all of the switches are equal to the
identity, that is, $\mathrm{id}_K \in K$.
It is worth noting, however, that the existence of a winning strategy does not
depend on a particular choice of the winning state.
Instead, we will see that a winning strategy is equivalent to a choice of moves
that will walk over all of the possible configuration states, regardless of the
choice of the adversary's spin.

\subsection{Surjective strategy}

We will begin by formalizing the notation of a winning strategy in a
generalized spinning switches puzzle. Informally, this is
a sequence of moves that the puzzle-solver can make that will put the switches
into every possible state, which ensures that the winning state is reached
regardless of the initial (hidden) state of the switches.

\begin{definition}
  A \textbf{surjective strategy} for $G \wr H$ is a finite sequence of elements
  in the base $K$,
  $\{k_i \in K\}_{i=1}^N$,
  such that for every sequence of elements in $H$, ${\{h_i \in H\}_{i=1}^N}$,
  \begin{equation}
    p(\{
      \underbrace{\mathrm{id}_{G \wr H}}_{m_0},
      \underbrace{(k_1, h_1)}_{m_1},
      \underbrace{(k_1, h_1)\cdot(k_2, h_2)}_{m_2},
      \dots,
      \underbrace{(k_1, h_1)\cdot(k_2, h_2)\cdots(k_N, h_N)}_{m_N}
    \}) = K.
  \end{equation}
  where $p \colon G \wr H \rightarrow K$ is the projection map from the
  wreath product onto its base.
\label{def:surjectiveStrategy}
\end{definition}

This definition is useful because it puts the problem into purely algebraic
terms. It is also useful because it abstracts away the initial state of the
switches: regardless of the initial state $k \in K$, the existence of a
surjective strategy means that its inverse $k^{-1} \in K$ appears in the sequence.
(This follows the convention that $\mathrm{id}_K \in K$ is designated as the
winning state. If $k'$ is chosen to be the winning state, then the sequence
must contain $k^{-1}k'$ .)

\begin{proposition}
  A finite sequence of moves is guaranteed to reach the winning
  state if and only if it is a surjective strategy.
\end{proposition}
\begin{proof}
  Without loss of generality, we will say that the winning state for the switches is
  $\mathrm{id}_K$.

  We will begin by assuming that $\{k_i \in K\}_{i=1}^N$ sequence of moves is
  guaranteed to reach the winning state, $\mathrm{id}_K$.
  In the puzzle, we have an initial (hidden) state, $k \neq \mathrm{id}_K$.
  Therefore, after the $i$-th move,
  the wreath product element that represents the state of the switches is
  \begin{equation}
    p\left((k, \mathrm{id}_H)\cdot(k_1, h_1)\cdot(k_2, h_2)\cdots(k_i, h_i)\right)
    = k \cdot p\left((k_1, h_1)\cdot(k_2, h_2)\cdots(k_i, h_i)\right),
  \end{equation}
  where the equality is due to associativity of the wreath product.
  We can factor out the first term because the ``spin'' is $\mathrm{id}_H$,
  which acts trivially:
  ${(k, \mathrm{id}_H) \cdot (k', h') = (kk', h')}$.

  Say that $j$ is the index at which the puzzle-solver gets the switches into
  the winning state. Then \begin{align}
    k \cdot p\left((k_1, h_1)\cdot(k_2, h_2)\cdots(k_j, h_j)\right) &= \mathrm{id}_K
    \\
            p\left((k_1, h_1)\cdot(k_2, h_2)\cdots(k_j, h_j)\right) &= k^{-1}.
    \label{eq:winningState}
  \end{align}
  Since $k \neq$ is arbitrary, there must exist such a $j$ for every $k \in K$
  and adversarial sequence $\{h_i \in H\}_{i=1}^N$. This means that the
  projection of the sequence of partial products of $\{k_i \in K\}_{i=1}^N$ is
  a surjection onto $K$ and therefore is a surjective strategy.

  Conversely, if $K$ is a surjective strategy, then for any initial state $k$ and sequence of
  adversarial moves $\{h_i \in H\}_{i=1}^N$, there exists some $j$ that
  satisfies Equation \eqref{eq:winningState}, and therefore reaches a winning
  state.
\end{proof}

It is also worth noting that this model can be thought of as a random model or an
adversarial model: the sequence $\{h_i \in H\}$ can be chosen
randomly at any point, or it can be chosen deterministically
after the sequence $\{k_i \in K\}$ is specified.

\subsection{Bounds on the length of surjective strategies}

One useful consequence of this definition is that it is quite straightforward to
prove certain propositions. For example, the minimum length for a surjective strategy
has a simple lower bound.
\begin{proposition}
  Every surjective strategy $\{k_i \in K\}_{i=1}^{N}$ is a sequence of
  length at least ${|K| - 1}$.
\end{proposition}
\begin{proof}
  This follows from an application of the Pigeonhole Principle, because the set
  \begin{equation}
    \{\mathrm{id}_{G \wr H}, (k_1, h_1), (k_1, h_1)\cdot(k_2, h_2), \ldots, (k_1, h_1)\cdot(k_2, h_2)\cdots(k_N, h_N)\}
  \end{equation}
  has at most $N+1$ elements. In order for the projection to be equal to $K$,
  \begin{equation}
    p(\{\mathrm{id}_{G \wr H}, (k_1, h_1), (k_1, h_1)\cdot(k_2, h_2), \cdots, (k_1, h_1)\cdot(k_2, h_2)\cdots(k_N, h_N)\}) = K,
  \end{equation} it must be the case that $N+1 \geq |K|$.
  Therefore $N \geq |K| - 1$.
\end{proof}

Minimal length surjective strategies are common, so we give them a name.
\begin{definition}
  A \textbf{minimal surjective strategy} for $G \wr H$ is a surjective strategy
  of length $N = |K| - 1.$
\end{definition}

In practice, every wreath product known to the author to have a surjective
strategy also has a known minimal surjective strategy.
In Section \ref{sec:OpenQuestions}, we ask whether this property always holds.





\section{Reductions}
\label{sec:Reductions}
In this section, we use three techniques to develop examples of generalized
spinning switches puzzles
that do not have surjective strategies:
directly, by a reduction on switches, or by a reduction on spinning.
\subsection{Puzzles known to have no surjective strategies}
Our richest collection of known puzzles without surjective strategies comes from
a theorem of Rabinovich, which models switches as a vector space over
a finite field.
\begin{theorem} \cite{Rabinovich2022}
  Assume that a finite ``spinning'' group $H$ acts linearly and faithfully on
  a collection of switches that behave like a vector space
  $V$ over a finite field $\mathbb F_q$ of characteristic $p$.
  Then the resulting puzzle has a surjective strategy if and only if $H$ is a
  $p$-group.
  \label{thm:Rabinovich}
\end{theorem}
It is worth noting that Rabinovich's switches are less general than arbitrary
finite groups, but the ``spinning'' is more general: in addition to permuting
the switches, the group action might add linear combinations of them as well.
\begin{example}
  By the theorem of Rabinovich \cite{Rabinovich2022},
  the game $\mathbb Z_2 \wr C_3$ does not have a
  surjective strategy. In Rabinovich's notation, the vector space of switches
  is $\mathbb Z_2^3$ over the field $\mathbb F_2 = \mathbb Z_2$, and
  $C_3$ has $3$ elements and therefore is not a $2$-group.
  \label{ex:NoSolutionZ2C3}
\end{example}

The wreath product $\mathbb Z_2 \wr C_3$ is perhaps the simplest example of a
generalized spinning switches puzzle without a surjective strategy,
so we will continue to use it as a basis of future examples.

\subsection{Reductions on switches}
With Theorem \ref{thm:Rabinovich} providing a family of wreath products without
surjective strategies, we now introduce a theorem that allows us to
describe large families of wreath products that also do not have surjective
strategies.
\begin{theorem}
  If $G \wr H$ does not have a surjective strategy and there exists a group $G'$
  and a surjective homomorphism $\varphi \colon G' \rightarrow G$,
  then ${G'} \wr H$ does not have a surjective
  strategy.
  \label{thm:SwitchReduction}
\end{theorem}
\begin{proof}
  We will prove the contrapositive, and suppose that $G' \wr H$ has base $K'$
  and a surjective strategy $\{k'_i \in K'\}_{i=1}^N$.

  The homomorphism
  $\varphi\colon G' \rightarrow G$
  extends coordinatewise to
  $\varphi \colon K' \rightarrow K$,
  which further extends in the first coordinate to $G \wr H$:
  $\varphi(k,h) := (\varphi(k), h)$.

  It is necessary to verify that $\varphi\colon G' \wr H \rightarrow G \wr H$
  is indeed a homomorphism.
  \begin{align*}
    \varphi((k'_\alpha, h_\alpha)) \cdot \varphi((k'_\beta, h_\beta))
    &= (\varphi(k'_\alpha), h_\alpha) \cdot (\varphi(k'_\beta), h_\beta) \\
    &= (\varphi(k'_\alpha)(h_\alpha\cdot\varphi(k'_\beta)), h_\alpha h_\beta) \\
    &= (\varphi(k'_\alpha)\varphi(h_\alpha\cdot k'_\beta), h_\alpha h_\beta) \\
    &= (\varphi(k'_\alpha(h_\alpha\cdot k'_\beta)), h_\alpha h_\beta) \\
    &= \varphi((k'_\alpha(h_\alpha\cdot k'_\beta), h_\alpha h_\beta)) \\
    &= \varphi((k'_\alpha, h_\alpha) \cdot (k'_\beta, h_\beta))
  \end{align*}

  Therefore the sequence $\{\varphi(k'_i) \in K\}_{i=1}^N$ is a
  surjective strategy on $G \wr H$, because the quotient map
  $\varphi \colon G' \rightarrow G$
  (and thus $\varphi \colon K' \rightarrow K$)
  is injective.
\end{proof}
\begin{example}
  We know that $\mathbb Z_2 \wr C_3$ does not have a surjective strategy.
  This means that $\mathbb Z_6 \wr C_3$ does not have a surjective strategy either,
  as illustrated in Figure \ref{fig:Z2C3}.
  \begin{figure}
    \center
    \includegraphics[width=\textwidth]{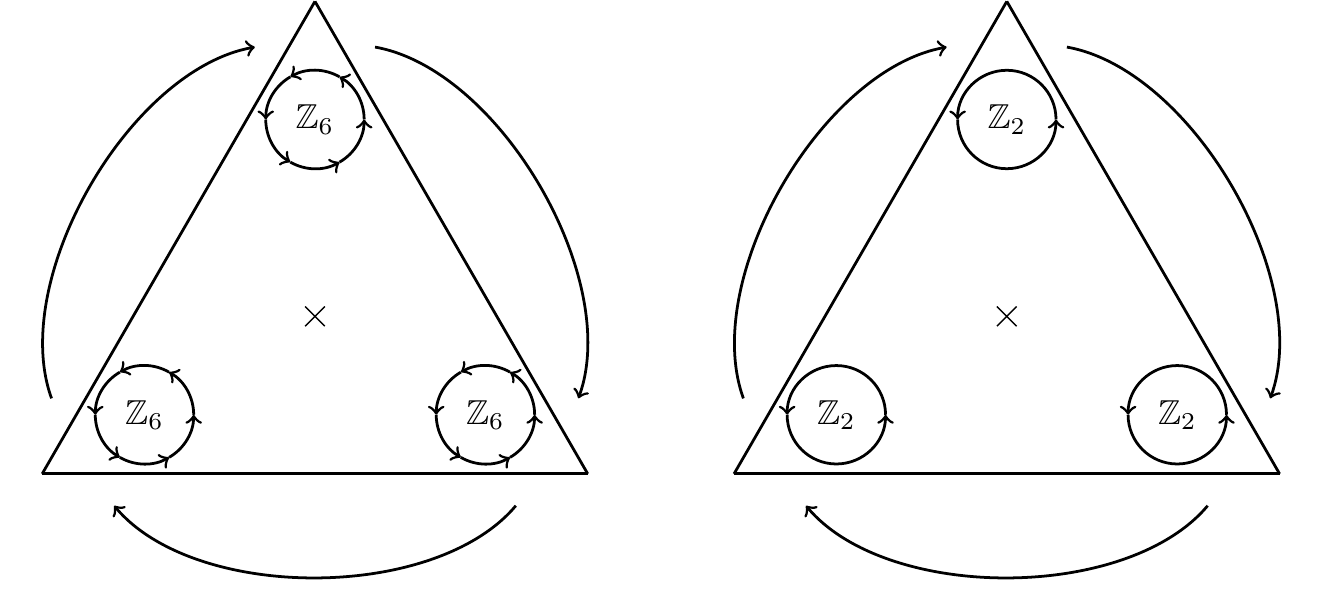}
    \caption[A reduction from a six-way switch to a two-way switch]{
      A reduction on switches:
      $\mathbb Z_6 \wr C_3$ reduces to $\mathbb Z_2 \wr C_3$,
      which is known not to have a surjective strategy.
    }
    \label{fig:Z2C3}
  \end{figure}
\end{example}
\subsection{Reductions on spinning}
We can do two similar reductions on the ``spinning'' group of a wreath product.
These theorems say that if a given wreath product $G \wr H$ does not have a
surjective strategy, then a similar wreath product $G \wr H'$ with a
``more complicated'' spinning group $H'$ will not have a surjective strategy
either.
\begin{theorem}
  If $G \wr H$ does not have a surjective strategy,
  and $\varphi \colon H \hookrightarrow H'$ is an embedding of $H$ into $H'$,
  then $G \wr H'$ does not have a surjective strategy.
  \label{thm:SpinReduction}
\end{theorem}
\begin{proof}
  Again we will prove the contrapositive.
  Assume that $G \wr H'$ does have a surjective strategy,
  $\{k_i\}_{i=1}^N$. Then by definition, for any sequence
  $\{h'_i\}_{i=1}^N$, the projection of the sequence \begin{equation}
    p(\{(k_1, h'_1)\cdot(k_2, h'_2)\cdots(k_i, h'_i)\}_{i=1}^N)
    = \left\{p\left((k_1, h'_1)\cdot(k_2, h'_2)\cdots(k_i, h'_i)\right)\right\}_{i=1}^N
    = K,
  \end{equation} and in particular this is true when $h'_i$ is restricted to be
  in the subgroup $\operatorname{Im}(\varphi) \leq H'$.
  Thus a surjective strategy for $G \wr H'$ is also a valid surjective strategy
  for $G \wr H$.
\end{proof}

\begin{figure}
  \center
  \includegraphics[width=\textwidth]{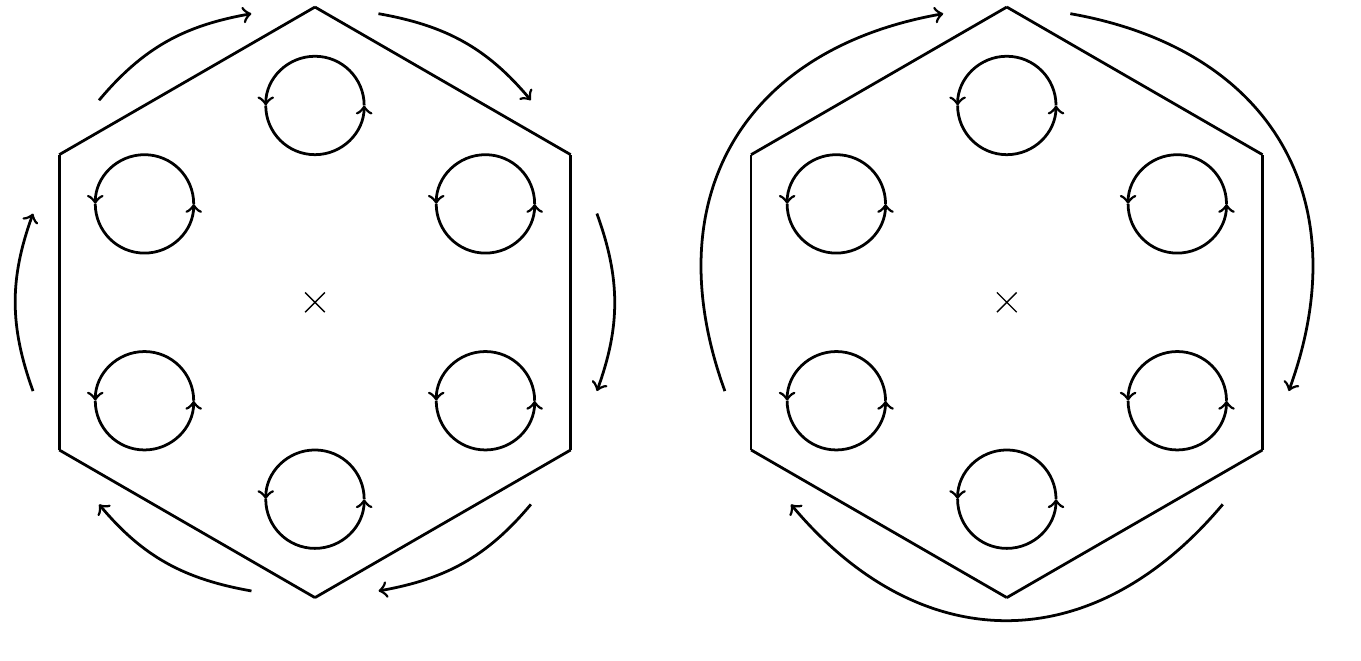}
  \caption[A reduction from $60^\circ$ rotations to $120^\circ$ rotations.]{
    If there were a solution to $\mathbb Z_2 \wr_{\Omega_6} C_6$,
  then there would be a solution to $\mathbb Z_2 \wr_{\Omega_6} C_3$.
  }
\label{fig:Z2C6}
\end{figure}

\begin{example}
  Consider the wreath product $\mathbb Z_2 \wr_{\Omega_6} C_6$ where
  $\Omega_6$ consists of six switches on the corners of a hexagon as
  illustrated in Figure \ref{fig:Z2C6}.
  We know that $\mathbb Z_2 \wr_{\Omega_6} C_3$ does not have a surjective
  strategy, by Theorem \ref{thm:Rabinovich}, therefore
  $\mathbb Z_2 \wr_{\Omega_6} C_6$ cannot have a surjective strategy either, by
  Theorem \ref{thm:SpinReduction}.
\end{example}

In the above example, we noted that $\mathbb Z_2 \wr_{\Omega_6} C_6$ does
not have a surjective strategy because ${\mathbb Z_2 \wr_{\Omega_6} C_3}$ is
known not to have one. However, there is another obstruction: the fact that
the generalized spinning switches puzzle $\mathbb Z_2 \wr C_3$ on a triangle
is known not to have a surjective strategy. That is, we cannot guarantee that
the upright triangle (which is formed by taking ``every other'' switch)
is ever in the all ``on'' state.

The following theorem abstracts this idea.
\begin{theorem}
  Suppose that $H$ and $H'$ are groups,
  $\varphi \colon H \hookrightarrow H'$ is an embedding of $H$ into $H'$, and
  let
  \begin{equation}
    \operatorname{Orb}(\omega)
      = \{\omega \cdot a : a \in \operatorname{Im}(\varphi) \} \subseteq \Omega
  \end{equation}
  be the (right) orbit of $\omega \in \Omega$ under $\operatorname{Im}(\varphi)$.
  If $G \wr_{\operatorname{Orb}(\omega)} H$ does not have a surjective strategy,
  then $G \wr_\Omega H'$ cannot have a surjective strategy either.
  \label{thm:SpinReduction2}
\end{theorem}
\begin{proof}
  We start by making the contrapositive assumption that $G \wr_\Omega H'$
  has a surjective strategy ${\{k_i \in K\}_{i=1}^N}$,
  and we consider the projection
  $p_\omega \colon K \rightarrow K_\omega$ where
  \begin{equation}
    K = \prod_{\omega' \in \Omega} G_{\omega'}
    \hspace{1cm}\text{and}\hspace{1cm}
    K_{\omega} = \prod_{\omega' \in \operatorname{Orb}(\omega)} G_{\omega'}.
  \end{equation}

  Then $\{p_\omega(k_i) \in K_\omega\}_{i=1}^N$ is a surjective strategy for
  $G \wr_{\operatorname{Orb}(\omega)} H$, since the projection is a surjective
  map.
\end{proof}

\begin{figure}
  \center
  \includegraphics[width=\textwidth]{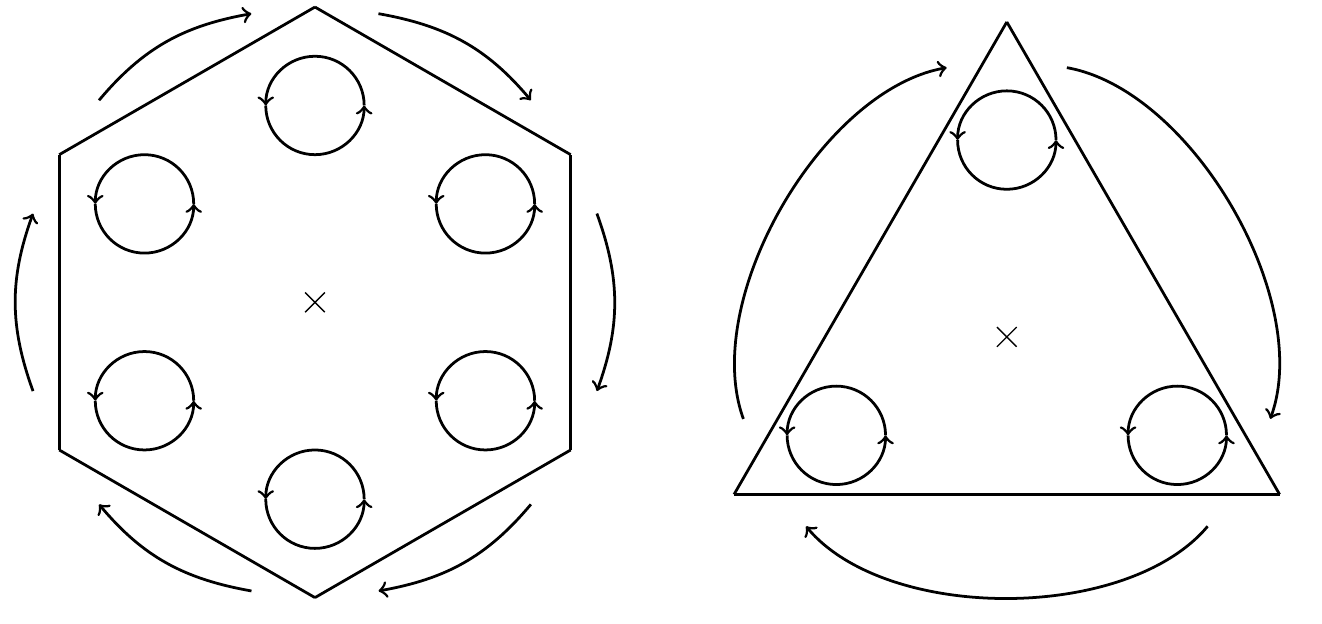}
  \caption[A reduction from a hexagonal table to a triangular table.]{
    We know that $\mathbb Z_2 \wr_\Omega C_6$ cannot have
    a surjective strategy, because that would imply a surjective strategy for
    $\mathbb Z_2 \wr_{\Omega'} C_3$, where $\Omega'$
    is the orbit of the top switch rotations of multiples of $120^\circ$.
  }
  \label{fig:Z2C6_2}
\end{figure}

\begin{example}
  We know that $\mathbb Z_2 \wr C_3$ does not have a surjective strategy.
  This means that $\mathbb Z_2 \wr C_6$ does not have a surjective strategy either,
  as illustrated in Figure \ref{fig:Z2C6_2}.
\end{example}

Now that we have proven that large families of wreath products do not have
surjective strategies, it is worthwhile to construct families of wreath products
that do have surjective strategies.
%
%
\section{Surjective strategies on \texorpdfstring{$p$}{p}-groups}
\label{sec:pGroupStrategy}
In this section, we will develop a broad family of surjective strategies,
namely those where $G$ and $H$ (and thus $G \wr H$) are $p$-groups.
\subsection{Surjective strategy decomposition}

Our first constructive theorem provides a technique that can be used to
construct surjective strategies for switches that behave like a group $G$ in
terms of a normal group and its corresponding quotient group.
\begin{theorem}
  The wreath product $G \wr H$ has a surjective strategy if there exists a
  normal subgroup $N \trianglelefteq G$ such that both $N \wr H$ and
  $G/N \wr H$ have surjective strategies.
\label{thm:surjectiveStrategyDecomposition}
\end{theorem}
\begin{proof}
  Let \[
    S_{G/N} = \{k_i^{G/N} \in K_{G/N}\}
    \hspace{1cm}\text{and}\hspace{1cm}
    S_{N} = \{k_i^N \in K_{N}\}
  \] denote the surjective strategies for ${G/N \wr H}$ and ${N \wr H}$
  respectively.

  We ultimately would like to interleave these two strategies,
  but $k_i^{G/N} \not\in K_G$. To find the appropriate analog,
  we partition $G$ into $[G : N] = m$ right cosets of $N$, \begin{equation}
    G = Ng_1 \sqcup Ng_2 \sqcup \dots \sqcup Ng_m,
  \end{equation} each with a chosen representative in $G$.
  Now we define a map $r \colon G/N \rightarrow G$ that chooses the chosen
  representative of the coset, and extends coordinatewise. We use this map to
  define a sequence $S = \{r(k_i^{G/N}) \in K_G\}$.

  We claim that these two sequences interleaved, $S_N \circledast S$, forms a
  surjective strategy for $G \wr H$. To prove this claim, we observe two facts:
  \begin{enumerate}
    \item Multiplying by elements of $S_N$ will not change cosets and will walk
    through every element of the current coset.
    \item Multiplying by elements of $S$ will walk through all cosets.
  \end{enumerate}
  Therefore the interleaved sequence will walk through all elements of each
  coset, and thus is surjective onto $K$.
\end{proof}

\subsection{Construction of surjective strategies on \texorpdfstring{$p$}{p}-groups}
With the decomposition from Theorem \ref{thm:surjectiveStrategyDecomposition}
established, we can now construct a surjective strategy on all finite $p$-groups.

\begin{theorem}
  If $H$ is a finite $p$-group that acts faithfully on $\Omega$,
  then the wreath product
  $G \wr H$ has a surjective strategy whenever $|G| = p^n$ for some $n$.
  \label{thm:pGroupStrategy}
\end{theorem}
\begin{proof}
  We will use the fact that if $|G| = p^n$,
  then either $G \cong \mathbb{Z}_p$ or $G$ is not simple.

  If $n = 1$, then $G \cong \mathbb{Z}_p \cong \mathbb F_p$, so there exists a
  surjective strategy by Theorem \ref{thm:Rabinovich}. This is because $H$
  permutes the coordinates of $V = \mathbb F_p^{|\Omega|} \cong K$, and so it is a
  linear action on the vector space, and so it has a surjective strategy by
  Theorem \ref{thm:Rabinovich}.

  Otherwise, $G$ is not simple. This means that $G$ has a proper normal subgroup
  $N$ of order $|N| = p^t$ (with $0 < t < n$)
  and a quotient $G/N$ with order $|G/N| = p^{n-t}$.
  Therefore, by induction on the exponent,
  whenever $G$ and $H$ (and thus $G \wr H$) are $p$-groups,
  $G \wr H$ has a surjective strategy.
\end{proof}

This resolves the situation for switches that behave like $p$-groups on
spinning groups that act faithfully and behave like $p$-groups. This allows us
to fully classify the case of abelian switches.

\subsection{A classification of puzzles with abelian switches}

\begin{theorem}
  If $G$ is an abelian group, and $H$ acts faithfully on $\Omega$, then
  $G \wr_\Omega H$ has a surjective strategy if and only if $G$ and $H$ are
  both $p$-groups for the same prime $p$.
\label{thm:classifyAbelianSwitches}
\end{theorem}
\begin{proof}
  One direction is clear: if $G$ and $H$ are both $p$-groups for the same prime
  $p$, then $G \wr H$ has a surjective strategy by Theorem \ref{thm:pGroupStrategy}.

  For the other direction, we will consider the contrapositive assumption,
  split into three cases, each time assuming that
  $G$ and $H$ are not both $p$-groups for the same prime $p$.

  First, assume that $G$ is a $p$-group, but $H$ is not a $p$-group. Then by
  Sylow's first theorem, there exists a subgroup $\hat{H} \leq H$ that is a $q$-group
  for some prime $q \neq p$. By Theorem \ref{thm:Rabinovich}, since $q \neq p$,
  $G \wr \hat{H}$ does not have a surjective strategy, thus it follows from
  the reduction on spinning given in Theorem \ref{thm:SpinReduction} that
  $G \wr H$ does not have a surjective strategy either.

  Second, assume that $H$ is a $p$-group, but $G$ is not a $p$-group. Similarly,
  by Sylow's first theorem, there exists a subgroup $\hat{G} \leq G$ that is
  a $q$-group for some prime $q \neq p$. Since $G$ is abelian, $\hat{G}$ is a
  normal subgroup. By Theorem \ref{thm:Rabinovich}, since $q \neq p$,
  $\hat G \wr H$ does not have a surjective strategy, thus it follows from
  the reduction on switches given in Theorem \ref{thm:SwitchReduction} that
  $G \wr H$ does not have a surjective strategy either.

  Lastly, assume that there does not exist any $p$ such that $G$ is a $p$-group
  or $H$ is a $p$-group. Then $|G|$ and $|H|$ must have multiple prime divisors,
  so $G$ has a normal subgroup, $\overline G$, that is a $q$-group for some prime $q$, and
  $H$ has a normal subgroup $\overline H$ that is an $r$-group for some prime $r \neq q$.
  $\overline G \wr \overline H$ does not have a surjective strategy by Theorem \ref{thm:Rabinovich},
  so $G \wr \overline H$ does not have a surjective strategy by Theorem \ref{thm:SwitchReduction},
  so $G \wr H$ does not have a surjective strategy by Theorem \ref{thm:SpinReduction}.
\end{proof}

\subsection{A folklore conjecture}
Here we note a conjecture from folklore, which---if true---implies that we have
\textit{almost} solved the problem in its full generality.
\begin{conjecture}[Folklore]
  Almost all groups are $2$-groups.
\end{conjecture}

There are both computational and theoretical bases for this conjecture.
According to the On-Line Encyclopedia of Integer Sequences \cite{oeis},
there are $A000001(2^{10}) = 49487367289$ groups of order $2^{10}$ and there are
$A063756(2^{11}-1) = 49910536613$ groups of order less than $2^{11}$.
This means that more than $99.15\%$ of the groups of order less than $2^{11}$
are of order $2^{10}$.

If this conjecture is true, then most types of switches have surjective
strategies on most kinds of faithful finite group actions.
Of course, while most finite groups may be $2$-groups, most mathematicians
are more interested in the other finite groups.
This next section develops two families of examples of surjective strategies
where the switches do not behave like $p$-groups.

%
%
\section{Surjective strategies on other wreath products}
\label{sec:OtherSurjectiveStrategies}

Other authors have given surjective strategies for various configurations of
generalized spinning switches strategies, but in all of these examples,
the wreath products themselves are $p$-groups: that is,
${|G \wr_\Omega H| = |G|^{|\Omega|} \cdot |H|}$, where $H$ acts faithfully.
In this section we introduce two families of wreath products that have surjective
strategies, but are not $p$-groups.

\subsection{The trivial wreath product, \texorpdfstring{$G \wr \mathbf{1}$}{G wreath trivial group}}
The simplest---and least interesting---way to construct a wreath product with
a surjective strategy is to remove the adversary (or randomness) altogether,
by letting the spinning group be the trivial group $H = \mathbf{1}$.
Additionally we will consider the case where $|\Omega| = 1$, so there is only
one switch.
Because the adversary cannot ``spin'' the switches at all, the puzzle-solver has
perfect information the entire time.
We will see that whenever $G$ is a finite group,
$G \wr \mathbf{1} \cong G$ has not just one surjective strategy, but many.

\begin{proposition}
  The wreath product $G \wr \mathbf{1}$ has $(|G|-1)!$ minimal surjective
  strategies.
  \label{prop:countingTrivialSS}
\end{proposition}
\begin{proof}
  There are $(|G| - 1)!$ permutations of $G \setminus \{\mathrm{id}_G\}$, and
  we claim each one corresponds to a minimal surjective strategy. Namely,
  if $(k_1, k_2, \dots, k_{|G|-1})$ is such a permutation, then
  the sequence $\{k'_i\}_{i=1}^{|G|-1}$ where $k'_1 = k_1$ and
  $k'_i = k_{i-1}^{-1}k_i$ is a surjective strategy on $G \wr \mathbf{1}$.

  Then we claim by induction that
  $(k'_1, \mathrm{id})\cdot(k'_2, \mathrm{id})\cdots(k'_j, \mathrm{id}) = (k_j, \mathrm{id})$.
  By construction, the base case is true when $j = 1$. If the claim holds up to
  $j-1$, then \begin{equation}
    \underbrace{
      (k'_1, \mathrm{id})\cdot(k'_2, \mathrm{id})\cdots(k'_{j-1}, \mathrm{id})
    }_{(k_{j-1}, \mathrm{id})}
    (k'_j, \mathrm{id})
    = (k_{j-1}, \mathrm{id})(k_{j-1}^{-1}k_j, \mathrm{id})
    = (k_j, \mathrm{id}),
  \end{equation} as desired.
  Thus the projection of the partial products is
  \begin{align}
    &p(\{
      \underbrace{\mathrm{id}_{G \wr \mathbf{1}}}_{m_0},
      \underbrace{(k'_1, h_1)}_{m_1},
      \underbrace{(k'_1, h_1)\cdot(k'_2, h_2)}_{m_2},
      \cdots,
      \underbrace{(k'_1, h_1)\cdot(k'_2, h_2)\cdots(k'_{|G|-1}, h_{|G|-1})}_{m_{|G|-1}}
    \}) \\
    & \hspace{1cm} =
    p(\{
      \mathrm{id}_{G \wr \mathbf{1}},
      (k_1, \mathrm{id}),
      (k_2, \mathrm{id}),
      \cdots,
      (k'_{|G|-1}, \mathrm{id})
    \}) \\
    & \hspace{1cm} = \{\mathrm{id}_{G \wr \mathbf{1}}, k_1, k_2, \cdots,k_{|G|-1}\} = K,
  \end{align}
  where $\{k_1, k_2, \dots, k_{|G|-1}\}$ spans
  $G \setminus \{\mathrm{id}_G\} \cong K \setminus \{\mathrm{id}_K\}$ by
  assumption.
\end{proof}

While the trivial wreath product is a useful example to keep in mind for
generating counterexamples, we are generally more interested in the situation
where the adversary permutes the switches to create uncertainty for the
puzzle-solver.

\subsection{Two interchangeable groups generated by involutions}

In this section, we will construct a surjective strategy for the generalized
spinning switches puzzle $G \wr H$
that consists of two switches each behaving like
a group $G$, generated by involutions, together with a nontrivial spinning
action $H$.
(See, for example, Figure \ref{fig:S3Switch} which illustrates switches that
behave like $S_3 = \langle (12), (13)\rangle$)

This strategy relies on the fact that because each generator is its own
inverse, applying a generator to the first switch or the second switch
has no effect on their difference.

This surjective strategy has two parts.
The first part ensures that the two switches have every possible difference.
The second part shows that we can get the first switch
(with respect to the projection onto $K$) to take on every possible value
without changing the difference between the switches.
\begin{theorem}
  Suppose that $G$ is a finite group that can be generated by involutions.
  Then the generalized spinning switches puzzle $G \wr C_2$ consisting of two
  interchangeable copies of $G$ has a surjective strategy.
  \label{thm:involutionGeneratedGroups}
\end{theorem}
\begin{proof}
  We start by writing $G$ in terms of its generators:
  $G = \langle t_1, t_2, \dots, t_N \rangle$, where $t_i^{-1} = t_i$.
  Because this is the generating set, there exists a finite sequence of
  transpositions $(t_{i_1}, t_{i_2}, \dots, t_{i_M})$
  such that the partial products of the sequence generate $G$: \begin{equation}
    G = \{\mathrm{id}_G, t_{i_1}, t_{i_1}t_{i_2}, \dots, t_{i_1}t_{i_2}\cdots t_{i_M}\}.
  \end{equation}

  We develop the strategy in two parts. First, we provide a strategy
  $A = \{\alpha_i \in K\}$ such that for any adversarial sequence $\{h_i \in H\}$
  and element $g \in G$, there exists an $i \geq 0$
  such that the {$i$-th} partial product,
  $(g_{i,1}, g_{i,2}) = p((\mathrm{id}_K, \mathrm{id}_H)\cdot(\alpha_1, h_1)\cdot(\alpha_2, h_2)\dots(\alpha_i, h_i))$,
  has a difference of $g$. That is, $g_{i,1}g_{i,2}^{-1} = g$.

  To do this, we define $\alpha_j = (t_{i_j}, \mathrm{id}_G) \in K$,
  and notice that the difference of the coordinates is the same whether we
  add $\alpha_j$ or $(180^\circ)\cdot \alpha_j$ to an element
  $(g_1, g_2) \in K$: \begin{equation}
    g_1(g_2t_{i_j})^{-1} = g_1t_{i_j}^{-1}g_2^{-1} = (g_1t_{i_j})g_2^{-1}
  \end{equation}

  Because the partial products of $t_{i_j}$ cover $G$, there exists some
  $t_{i_1} t_{i_2}\dots t_{i_k} = g_1^{-1}gg_2$, so that
  \begin{equation}
    g_1(t_{i_1} t_{i_2}\dots t_{i_k})g_2^{-1}
    = g_1(g_1^{-1}gg_2)g_2^{-1}
    = g,
  \end{equation}
  as desired.

  Next, we give a strategy $B = \{\beta_i\}$ such that for any adversarial
  sequence $\{h_i \in H\}$ and element $g \in G$, there exists an $i \geq 0$
  such that the {$i$-th} partial product,
  \begin{equation}
    p((\mathrm{id}_K, \mathrm{id}_H)\cdot(\alpha_1, h_1)\cdot(\alpha_2, h_2)\cdots(\alpha_i, h_i))
  \end{equation}
  has a first coordinate equal to $g$.

  To do this, we define $\beta_j = (t_{i_j}, t_{i_j})$. This strategy is
  invariant up to actions of $H$, so we can see that regardless of the initial
  state $(g_1, g_2) \in K$, there exists some $k$ such that
  \begin{equation}
    \beta_1 \beta_2 \cdots \beta_k = (g_1^{-1}g, g_1^{-1}g)
  \end{equation} and therefore
   $(g_1, g_2) \beta_1 \beta_2 \cdots \beta_j = (g, g_2g_1^{-1}g)$, as desired.

   It is important to note that applying $\beta_j$ does not affect the
   difference: $(g_1, g_2) \beta_j = (g_1t_{i_j}, g_2t_{i_j})$ has a
   difference of
   $(g_1t_{i_j})(g_2t_{i_j})^{-1} = g_1t_{i_j}t_{i_j}^{-1}g_2^{-1} = g_1g_2^{-1}$.

  Now by interleaving these two strategies, we see that the partial products
  of $B \circledast A$ hit every possible first letter and every possible
  difference for every first letter,
  therefore the projection of the set partial products of
  $B \circledast A$ covers $K$ for all adversarial sequences, so
  $B \circledast A$ is a surjective strategy.
\end{proof}
We will illustrate this idea explicitly letting $G$ be the smallest nonabelian
group of composite order, the symmetric group on three letters $S_3$, which is
isomorphic to the dihedral group of the triangle, $D_6$.
\begin{example}
  Note that $S_3 = \langle(12), (13)\rangle$ is generated by involutions, and
  that the partial products of the sequence $((12),(13),(12),(13),(12))$
  cover $S_3$.

  As a bit of notation, for a permutation $\pi \in S_3$, let
  $\pi_1 = (\pi, \mathrm{id}_{S_3}) \in K$
  and
  $\pi_2 = (\pi, \pi) \in K$, corresponding to sequences $B$ and $A$
  respectively in the above theorem.
  Then the following is a (minimal) surjective strategy on $S_3 \wr C_2$:
  \begin{singlespace}
  \begin{align*}
    (12)_2(13)_2(12)_2(13)_2(12)_2 \\
    &(12)_1 \\
    (12)_2(13)_2(12)_2(13)_2(12)_2 \\
    &(13)_1 \\
    (12)_2(13)_2(12)_2(13)_2(12)_2 \\
    &(12)_1 \\
    (12)_2(13)_2(12)_2(13)_2(12)_2 \\
    &(13)_1 \\
    (12)_2(13)_2(12)_2(13)_2(12)_2 \\
    &(12)_1 \\
    (12)_2(13)_2(12)_2(13)_2(12)_2
  \end{align*}
  \end{singlespace}
  \label{ex:TwoSymmetricGroups}
\end{example}
It is natural to ask which
spinning groups $H$ correspond to to a generalized spinning switches puzzle
$S_n \wr H$ with a surjective strategy.
We can use a reduction on switches to rule out a large number of these.
\begin{proposition}
  For $n > 1$, $S_n \wr H$ does not have a surjective strategy whenever
  $|H|$ is not a power of $2$.
\end{proposition}
\begin{proof}
  The alternating group $A_n$ is an index $2$ subgroup of $S_n$, so $A_n$ is
  normal, and $S_n/A_n \cong \mathbb Z_2$.
  Since we know that $\mathbb Z_2 \wr H$ has no surjective strategy when
  $|H|$ is not a power of $2$,
  by the reduction in Theorem \ref{thm:SwitchReduction}, $S_n \wr H$ does
  not have a surjective strategy.
\end{proof}

Many groups are generated by transpositions, including $22$ of the $26$
sporadic simple groups and the alternating groups $A_5$ and $A_n$ for $n > 9$,
as shown by Mazurov and Nuzhin respectively.
\begin{theorem}\cite{Mazurov2003}
  Let $G$ be one of the 26 sporadic simple groups.
  The group $G$ cannot be generated by three involutions two of which commute
  if and only if $G$ is isomorphic to $M_{11}$, $M_{22}$, $M_{23}$, or $M^c_L$.
\end{theorem}
\begin{theorem}\cite{Nuzhin1992}
  The alternating group $A_n$ is generated by three involutions,
  two of which commute, if and only if $n \geq 9$ or $n = 5$.
\end{theorem}

Nuzhin provides other families of groups that are generated by three
involutions, two of which commute, in subsequent papers
\cite{Nuzhin0,Nuzhin1,Nuzhin2}.

Thus, we have characterized finite wreath products with surjective strategies
in many cases:
wreath products that are $p$-groups,
trivial wreath products, and
wreath products $G \wr C_2$ where $G$ is generated by involutions.
In the next section, we provide even more general constructions and ask
more specific questions.

%
%
\section{Generalizations and open questions}
\label{sec:OpenQuestions}
In this section, we provide conjectures and suggest open questions about
the structure of surjective strategies when they exist,
further generalizations of spinning switches puzzles,
and, lastly, introduce a notion of an infinite surjective strategy for infinite
wreath products.

The ultimate open question is a full classification
of finite wreath products with surjective strategies.
\begin{openquestion}
  What finite wreath products, $G \wr H$, have a surjective strategy?
\end{openquestion}
\subsection{Switches generated by elements of prime power order}

In Theorem \ref{thm:involutionGeneratedGroups},
we constructed a strategy for $G \wr C_2$,
when $G$ can be generated by elements of order $2$.
We conjecture that that there is a broader construction for the case where
the adversary can act on the switches with a group of order $2^\ell$.

\begin{conjecture}
  When $G$ is a finite group generated by involutions, and
  $H$ is $2$-group that acts faithfully on the set of switches,
  there exists a surjective strategy for $G \wr H$.
\end{conjecture}

We can also ask about three switches that are generated by elements of order
$3$ on the corners of a triangular table. In particular, the alternating group
is such a group, and we conjecture that it has a surjective strategy.

\begin{conjecture}
  There exists a surjective strategy for $A_n \wr C_3$.
\end{conjecture}

Putting these two conjectures together, we boldly predict a large family
of wreath products with surjective strategies.
\begin{conjecture}
  If $G$ can be generated by elements of order $p^n$, and $H$ is a $p$-group
  acting faithfully on the set of switches $\Omega$, then $G \wr_\Omega H$ has
  a surjective strategy.
\end{conjecture}

\subsection{Palindromic surjective strategies}
In all known examples, when there exists a surjective strategy $S$,
we also know of a \textit{palindromic} surjective strategy
$S' = \{k'_i \in K\}_{i=1}^N$
such that $k'_i = k'_{N-i+1}$ for all $i$.

\begin{conjecture}
  Whenever $G \wr H$ has a surjective strategy, it also has a palindromic surjective
  strategy.
\end{conjecture}

If this conjecture is false, we suspect a counterexample can be found in the
case of the trivial wreath product, $G \wr \mathrm{1} \equiv G$.

\subsection{Quasigroup switches}
\label{sub:quasigroupSwitches}
In Section \ref{sec:GeneralizingSwitches}, we argued for modeling switches as
finite groups because of some desirable properties: \begin{enumerate}
  \item Closure. Regardless of which state a switch is in, every move results
    in a valid state.
  \item Identity. We do not have to move a switch on a given turn.
  \item Inverses. If a switch is off, we can always turn it on.
\end{enumerate}

In the list, we also included the axiom of associativity for three reasons:
switches in practice typically have associativity,
groups are easier to model than quasigroups with identity, and
associativity makes defining a ``surjective strategy'' simpler.

However, one could design a non-associative switch and the puzzle would still
be coherent. This is because the process of
a generalized spinning switches puzzle is naturally
``left associative'' in the language of programming language theory:
we are always ``stacking'' our next move onto the right.
As such, it is worth noting the slightly more general way of modeling switches:
as quasigroups with identity, called \textit{loops}.

In particular, we are interested in the smallest loop that is not a group
\cite{MSELoop}, which we denote as
$\mathcal{L} = (\{1, a, b, c, d\}, \ast)$,
and describe via its multiplication table, a Latin square of order $5$.
(Notice that $(ab)d = a \neq d = a(bd)$.)
\begin{singlespace}
\begin{equation}
  \begin{array}{c|ccccc}
    \ast & 1 & a & b & c & d \\
    \hline
      1 & 1 & a & b & c & d \\
      a & a & 1 & c & d & b \\
      b & b & d & 1 & a & c \\
      c & c & b & d & 1 & a \\
      d & d & c & a & b & 1
  \end{array}
\end{equation}
\end{singlespace}

\begin{conjecture}
  There exists a nontrivial adversarial group $H$ such that the generalized
  spinning switches puzzle with switches that behave like $\mathcal{L}$ has
  a winning strategy for the puzzle-solver.
\end{conjecture}

\subsection{Expected number of turns}
In practice, a puzzle-solver can get the switches into the winning state by
playing randomly.
Random play will \textit{eventually} turn on the lightbulb with probability $1$,
due to the finite number of configurations and the law of large numbers.

Thus we drop the requirement of a finite strategy and
ask about the expected value of the number of turns given
various sequences of moves.
Notice that this is an interesting question even (perhaps especially)
in the context of generalized spinning switches puzzles that do not have a
surjective strategy.

Winkler \cite{Winkler2021}
notes in the solution ``Spinning Switches'':
\begin{quote}
  Although no fixed number of steps can guarantee turning the bulb on in the
  three-switch version [with two-way switches],
  a smart randomized algorithm can get the bulb on in at most $5 \frac{5}{7}$
  steps on average, against any strategy by an adversary who sets the initial
  configuration and turns the platform \cite{Winkler2021}.
\end{quote}

%

% ...
%
A basic model for computing the expected number of turns assumes that
the initial hidden state $k \in K$ is not the winning state $\mathrm{id}_K$,
and that the adversary's ``spins'' are independent and identically distributed
uniformly random elements $h_j \in H$.

\begin{proposition}
  If the puzzle-solver chooses $k_j \in K \setminus \{\mathrm{id}_K\}$ uniformly
  at random (that is, never choosing the ``do nothing'' move)
  then the distribution of the resulting state after each turn will be uniformly
  distributed
  among the $|K| - 1$ different states. The probability of the resulting state
  being the winning state is
  \begin{equation}
    \mathbb{P}(p((k_1, h_1)\dots(k_j, h_j))=k^{-1}\ \mid\ p((k_1, h_1)\dots(k_{j-1}, h_{j-1})\neq k^{-1})) = \frac{1}{|K| - 1},
  \end{equation} and the expected number of moves is $|K| - 1$.
\label{prop:randomStrategy}
\end{proposition}
\begin{proof}
  Because the new states are in $1$-to-$1$ correspondence with the elements of
  $K \setminus \{\mathrm{id}_K\}$, and since
  $k_j \in K \setminus \{\mathrm{id}_K\}$
  is chosen uniformly at random, the projection of the partial product
  ${p((k_1, h_1)(k_2, h_2)\cdots(k_j, h_j))}$
  is uniformly distributed among all elements of
  ${K \setminus \{p((k_1, h_1)(k_2, h_2)\cdots(k_{j-1}, h_{j-1}))\}}$.
  The expected value is $|K| - 1$ because the number of turns follows a
  geometric distribution with parameter ${(|K| - 1)^{-1}}$.
\end{proof}

When a generalized spinning switches puzzle has a \textit{minimal} surjective
strategy, we can \textit{guarantee} that we turn on the light within $|K|-1$
moves, so we can certainly solve in fewer than $|K|-1$ moves on average.

\begin{proposition}
If the generalized spinning switches puzzle, $G \wr H$, has a minimal surjective
strategy, then the expected number of moves is $|K|/2$.
\end{proposition}
\begin{proof}
  If there exists a surjective strategy of length $|K| - 1$,
  then for each adversarial strategy ${\{h_i \in H\}_{i=1}^{|K| - 1}}$
  the projection of the sequence of partial products of moves induces
  a permutation of $K \setminus \{\mathrm{id}_K\}$.

  If the initial hidden state $k$ is chosen uniformly at random, then
  $k^{-1}$ is equally likely to occur at any position in this permutation,
  so the index of the winning state is uniform on
  ${\{1, 2, \dots, |K| - 1\}}$ and the expected number of moves is
  \begin{equation}
    \frac{1 + 2 + \dots + |K| - 1}{|K| - 1} = |K|/2.
  \end{equation}
\end{proof}

As this proposition suggests, we can always come up with a strategy that does
better than the uniformly random strategy in Proposition \ref{prop:randomStrategy}.

\begin{proposition}
  For every generalized spinning switches puzzle, $G \wr H$ such that
  $|K| > 2$, there always exists a (perhaps infinite) sequence
  whose expected number of moves is strictly less than $|K| - 1$.
\end{proposition}
\begin{proof}
  When $|K| > 2$, we can always improve on the random strategy in
  Proposition \ref{prop:randomStrategy}, by avoiding the move of
  $(g,\cdots, g) \in K$ followed by $(g^{-1}, \cdots  g^{-1}) \in K$,
  because the second move will put us into a previous state and so will
  turn on the lightbulb with probability $0$.
\end{proof}

We conjecture that this technique can be extended, and that the puzzle-solver
can always do asymptotically better than randomly guessing.

\begin{conjecture}
  There exists a constant $\frac{1}{2} < c < 1$ such that for all
  (finite) wreath products $G \wr H$ with sufficiently large $|K|$,
  the expected number of moves is less than $c|K|$.
\end{conjecture}

\subsection{Shortest surjective strategies}
Based on all of the examples that we know of, we conjecture that there exist
minimal surjective strategies whenever there exists a surjective strategy of any
length.
\begin{conjecture}
  Whenever $G \wr H$ has a surjective strategy, it also has a minimal surjective
  strategy $\{k_i \in K\}_{i=1}^{|K| - 1}$.
\end{conjecture}

On the other extreme, we have a weaker conjecture: whenever a wreath product has
a surjective strategy, we can provide an upper bound for its shortest surjective
strategy.

\begin{conjecture}
  Let $K/H$ be the set of equivalence classes of $K$ up to the action of $H$.
  Then if $G \wr H$ has a surjective strategy, it always has a surjective strategy
  of length $N < 2^{|K/H|-1}$.
\end{conjecture}

\subsection{Counting surjective strategies}
The counting problem analog to the decision problem
``does $G \wr H$ have a surjective strategy'' is of obvious interest to the
author. It is interesting to count both the number of surjective
strategies of length $N$, and the number of such surjective strategies
\textit{up to the action of} $H$. We might also be interested in the number of
palindromic surjective strategies, or the number of surjective strategies
satisfying another desirable criterion.

In the case of the trivial wreath product, $G \wr \mathbf{1}$, we saw in
Proposition \ref{prop:countingTrivialSS} that there are ${(|G| - 1)!}$
surjective strategies. (And since the group action is trivial, there are also
this many strategies up to group action.)

\begin{openquestion}
  Given a wreath product $G \wr H$ how many surjective strategies of length $N$
  does it have? How many up to the action of $H$? How many are palindromic?
\end{openquestion}

\begin{proposition}
  The trivial wreath product $S_3 \wr \bf{1}$ has $12$ palindromic surjective
  strategies of length $5$:
  \captionsetup{type=table}
  \begin{singlespace}
  \begin{alignat*}{5}
    (1~2), && (1~3), && (1~2), && (1~3), && (1~2). \\
    (1~2), && (2~3), && (1~2), && (2~3), && (1~2). \\
    (1~3), && (1~2), && (1~3), && (1~2), && (1~3). \\
    (1~3), && (2~3), && (1~3), && (2~3), && (1~3). \\
    (1~2~3), && \hspace{1em} (1~2~3), && \hspace{1em} (1~2), && \hspace{1em} (1~2~3), && \hspace{1em} (1~2~3). \\
    (1~2~3), && (1~2~3), && (1~3), && (1~2~3), && (1~2~3). \\
    (1~2~3), && (1~2~3), && (2~3), && (1~2~3), && (1~2~3). \\
    (1~3~2), && (1~3~2), && (1~2), && (1~3~2), && (1~3~2). \\
    (1~3~2), && (1~3~2), && (1~3), && (1~3~2), && (1~3~2). \\
    (1~3~2), && (1~3~2), && (2~3), && (1~3~2), && (1~3~2). \\
    (2~3), && (1~2), && (2~3), && (1~2), && (2~3). \\
    (2~3), && (1~3), && (2~3), && (1~3), && (2~3).
  \end{alignat*}
  \end{singlespace}
  \captionof{table}[Palindromic walks on $S_3$.]{
    The 12 palindromic surjective strategies of length $5$ on $S_3 \wr \mathbf{1}$.
  }
\end{proposition}
\begin{proof}
  The search space is small here, so this was computed naively by brute force.
\end{proof}
\subsection{Restricted spinning}
Another way that we could generalize a spinning switches puzzle
is by restricting the adversary's moves.
For instance, we could modify the puzzle in such a way that the adversary can
only spin the switches every $k$ turns.
For every finite setup $G \wr H$, there exists $k \in \mathbb N$ such that the
puzzle-solver can win.
(For example, we can always take $k > |K|$ so that the puzzle-solver can
do a walk of $K$.)

\begin{openquestion}
  How does one compute the minimum $k$ such that the puzzle solver has a
  surjective strategy of $G \wr H$ given that the adversary's sequence
  $\{h_i \in H\}$ is constrained so that $h_i = \mathrm{id}_H$ whenever ${i \not\equiv 0 \bmod k}$?
\end{openquestion}

\subsection{Multiple winning states}
One assumption that we made about our switches is that they have a single
``on'' state. Of course, we might conceive of a switch that has $m$ possible
states and $k \leq m$ ``on'' states.

\begin{openquestion}
  Given a generalized spinning switches puzzle $G \wr H$ where each switch has
  a set of ``on'' states $\mathcal{O}_G \subseteq G$, when is it possible for
  the puzzle-solver to have a finite strategy that guarantees the switches will
  get to a winning state?
\end{openquestion}

\subsection{Nonhomogeneous switches}

Another way to generalize the spinning switches puzzle is by allowing different
sorts of switches together in the same puzzle.
For instance, we could imagine a square board
containing four buttons: two 2-way switches, a 3-way switch, and a 5-way switch.
Can a puzzle like this be solved?
It is important that all of these buttons have the same ``shape'', that is
there is some group $G$ with a surjective homomorphism onto each of them.



One way to formalize this generalization is as follows:
\begin{definition}
  A \textbf{nonhomogeneous generalized spinning switches puzzle}
  consists of a triple \begin{itemize}
    \item A wreath product \(
      \mathbb F_k \wr_\Omega H
    \) of the free group on $k$ generators by a ``rotation'' group
    with a base denoted \(K = \prod_{\omega \in \Omega} F_{k, \omega}\),
    \item a product of finite groups denoted
    $\hat G = \prod_{\omega \in \Omega} G_\omega$ where each group is specified
    by a presentation with $k$ generators and any number of relations:
    \({
      G_\omega = \langle g^\omega_1, g^\omega_2, \dots, g^\omega_k\ |\ R_\omega\rangle,
    }\) and
    \item a corresponding sequence of evaluation maps
    $e_\omega \colon \mathbb F_{k,\omega} \rightarrow G_\omega$,
    that send generators in $F_{k,\omega}$ to the corresponding generators in
    $G_\omega$. This can be induced coordinatewise to a map
    $e_\Omega \colon K \rightarrow \hat{G}$.
  \end{itemize}
\end{definition}

When all of the copies of $G_\omega$ are isomorphic, this essentially simplifies
to the original definition.

The analogous definition of a surjective strategy becomes more complicated.

\begin{definition}
  Let $(\mathbb F_k \wr_\Omega H, \hat{G}, e_\Omega)$ be a nonhomogeneous
  generalized spinning switches puzzle.

  Then a \textbf{nonhomogeneous surjective strategy} is a sequence
  $\{k_i \in K\}_{i=1}^{N}$
  such that for each adversarial sequence ${\{h_i \in H\}_{i=1}^{N}}$
  the induced pro\-jec\-tion/eval\-u\-a\-tion map
  $e_\Omega \circ p \colon \mathbb F_k \wr_\Omega H \rightarrow \hat{G}$
  on the partial products of $\{(k_i, h_i) \in \mathbb F_k\}$
  covers $\hat{G}$.
\end{definition}

\begin{proposition}
  In the specific case that $\Omega = [n]$, $H \subseteq S_n$, and
  $\hat{G} = \prod_{\omega \in \Omega} G_\omega$ is a product of cyclic groups of
  pairwise coprime order,
  the nonhomogeneous generalized spinning switches puzzle has a surjective
  strategy, namely $\{(1,1,\dots,1) \in K\}_{i = 1}^{|K'| - 1}$.
\end{proposition}

\begin{proof}
  By the fundamental theorem of abelian groups, $\hat{G}$ is cyclic and is
  generated by \begin{equation}
    \hat{G} = \langle(1_{C_{k_1}}, 1_{C_{k_2}}, \dots, 1_{C_{k_{|\Omega|}}})\rangle.
  \end{equation}
  Thus, the sequence $\{(1,1,\dots,1) \in K\}_{i=1}^{|K'|-1}$ is a surjective
  strategy because it is a fixed point under $H$, and its image is a generator
  of $\hat{G}$.
\end{proof}

\begin{openquestion}
  Which nonhomogeneous generalized spinning switches puzzles have
  a nonhomogeneous surjective strategy?
\end{openquestion}

\subsection{Infinite surjective strategies}
When we first introduced the notion of a surjective strategy in
Definition \ref{def:surjectiveStrategy},
we defined it to be a finite sequence on finite wreath products.
However, we can expand the definition to (countably) infinite wreath products
by allowing for infinite sequences of moves in $K$.
In particular, we can extend this definition to settings where
switches have a countably infinite number of states,
where there is a countably infinite number of switches,
or both.
To keep $K$ countable in the latter cases,
we use the restricted wreath product, where
$K \cong \bigoplus_{\omega \in \Omega} G_\omega$ is defined to be a direct
sum instead of a direct product.

\begin{definition}
  A \textbf{infinite surjective strategy} on an infinite wreath product $G \wr H$
  is a sequence $\{k_i \in K\}_{i=1}^\infty$ such that for all $k \in K$ and
  all infinite sequences ${\{h_i \in H\}_{i=1}^\infty}$,
  there exists some $N \geq 0$ such that the projection \begin{equation}
    p((k_1, h_1)\cdot(k_2, h_2)\cdots(k_N, h_N)) = k^{-1}.
  \end{equation}
\end{definition}

We claim, but do not prove that $G \wr_\Omega C_2$ has an infinite surjective
strategy in the following three settings: \begin{enumerate}
  \item $|\Omega| = 2$ and $G \cong (\mathbb N_{\geq 0}, \wedge)$ where $\wedge$ is the bitwise XOR operator.
  \item $\Omega \cong \mathbb N_{\geq 0}$ as a set, and $K = \bigoplus_{i=0}^\infty \mathbb Z_2^{(i)}$ where $C_2^{(i)} \cong \mathbb Z_2$.
  \item $\Omega \cong \mathbb N_{\geq 0}$ as a set, and $K = \bigoplus_{i=0}^\infty (\mathbb N_{\geq 0}, \wedge)$ for all $i$.
\end{enumerate}

\printbibliography
\end{document}